\documentclass[11pt,reqno,oneside]{amsproc}
\title[Kuramoto-Sivashinsky $L^{5/6}$ bound]{An $L^{5/6}$ bound for the periodic Kuramoto-Sivashinsky equation}

\author[B.~Ingimarson]{Benjamin Ingimarson}
\address{Department of Mathematics, University of Southern California, Los Angeles, CA 90089}
\email{ingimars@usc.edu}

\author[I.~Kukavica]{Igor Kukavica}
\address{Department of Mathematics, University of Southern California, Los Angeles, CA 90089}
\email{kukavica@usc.edu}

\author[Q.~Xu]{Qi Xu}
\address{Department of Mathematics, University of Southern California, Los Angeles, CA 90089}
\email{xuqi@usc.edu}

  \chardef\forshowkeys=0
  \chardef\showllabel=0
  \chardef\refcheck=0
  \chardef\sketches=0
  \chardef\figures=1

\usepackage{enumitem,mathtools}
\usepackage{datetime}
\usepackage{fancyhdr}
\usepackage{comment}
\usepackage{listings}

\allowdisplaybreaks
\ifnum\forshowkeys=1

  \usepackage[notref,notcite,color]{showkeys}
\fi
\usepackage[margin=1in]{geometry}
\usepackage{amsmath, amsthm, amssymb}
\usepackage{times}
\usepackage{graphicx}
\usepackage[usenames,dvipsnames,svgnames,table]{xcolor}
\usepackage{marginnote}
\usepackage[unicode,breaklinks=true,colorlinks=true,linkcolor=blue,urlcolor=blue,citecolor=blue]{hyperref}
\usepackage[most]{tcolorbox}
\usepackage{tikz}

\ifnum\refcheck=1
  \usepackage{refcheck}
\fi

\usepackage{import}
\usepackage{xifthen}
\usepackage{pdfpages}
\usepackage{transparent}

\begin{document}
\def\YY{X}
\def\OO{\mathcal O}
\def\SS{\mathbb S}
\def\CC{\mathbb C}
\def\RR{\mathbb R}
\def\TT{\mathbb T}
\def\TTL{\mathbb T}
\def\ZZ{\mathbb Z}
\def\HH{\mathbb H}
\def\RSZ{\mathcal R}
\def\LL{\mathcal L}
\def\SL{\LL^1}
\def\ZL{\LL^\infty}
\def\GG{\mathcal G}
\def\tt{\langle t\rangle}
\def\erf{\mathrm{Erf}}
\def\mgt#1{\textcolor{magenta}{#1}}
\def\ff{\rho}
\def\gg{G}
\def\sqrtnu{\sqrt{\nu}}
\def\ww{w}
\def\ft#1{#1_\xi}
\def\les{\lesssim}
\def\lec{\lesssim}
\def\ges{\gtrsim}
\def\gec{\gtrsim}
\renewcommand*{\Re}{\ensuremath{\mathrm{{\mathbb R}e\,}}}
\renewcommand*{\Im}{\ensuremath{\mathrm{{\mathbb I}m\,}}}
\ifnum\showllabel=1
 \def\llabel#1{\marginnote{\color{lightgray}\rm\small(#1)}[-0.0cm]\notag}
\else
 \def\llabel#1{\notag}
\fi
\newcommand{\norm}[1]{\left\|#1\right\|}
\newcommand{\nnorm}[1]{\lVert #1\rVert}
\newcommand{\abs}[1]{\left|#1\right|}
\newcommand{\NORM}[1]{|\!|\!| #1|\!|\!|}

\newtheorem{Theorem}{Theorem}[section]
\newtheorem{Corollary}[Theorem]{Corollary}
\newtheorem{Proposition}[Theorem]{Proposition}
\newtheorem{Lemma}[Theorem]{Lemma}
\newtheorem{Remark}[Theorem]{Remark}
\newtheorem{definition}{Definition}[section]

\def\theequation{\thesection.\arabic{equation}}
\numberwithin{equation}{section}
\definecolor{mygray}{rgb}{.6,.6,.6}
\definecolor{myblue}{rgb}{9, 0, 1}
\definecolor{colorforkeys}{rgb}{1.0,0.0,0.0}
\newlength\mytemplen
\newsavebox\mytempbox
\makeatletter
\newcommand\mybluebox{%
    \@ifnextchar[
       {\@mybluebox}%
       {\@mybluebox[0pt]}}
\def\@mybluebox[#1]{%
    \@ifnextchar[
       {\@@mybluebox[#1]}%
       {\@@mybluebox[#1][0pt]}}
\def\@@mybluebox[#1][#2]#3{
    \sbox\mytempbox{#3}%
    \mytemplen\ht\mytempbox
    \advance\mytemplen #1\relax
    \ht\mytempbox\mytemplen
    \mytemplen\dp\mytempbox
    \advance\mytemplen #2\relax
    \dp\mytempbox\mytemplen
    \colorbox{myblue}{\hspace{1em}\usebox{\mytempbox}\hspace{1em}}}
\makeatother

\def\pv{\text{p.v.}}
\def\bnew{\color{red}}
\def\enew{\color{black}}
\def\bold{\color{blue}}
\def\eold{\color{black}}
\def\rr{r}
\def\weaks{\text{\,\,\,\,\,\,weakly-* in }}
\def\inn{\text{\,\,\,\,\,\,in }}
\def\cof{\mathop{\rm cof\,}\nolimits}
\def\Dn{\frac{\partial}{\partial N}}
\def\Dnn#1{\frac{\partial #1}{\partial N}}
\def\tdb{\tilde{b}}
\def\tda{b}
\def\qqq{u}
\def\lat{\Delta_2}
\def\biglinem{\vskip0.5truecm\par==========================\par\vskip0.5truecm}
\def\inon#1{\hbox{\ \ \ \ \ \ \ }\hbox{#1}}                
\def\onon#1{\inon{on~$#1$}}
\def\inin#1{\inon{in~$#1$}}
\def\FF{F}
\def\andand{\text{\indeq and\indeq}}
\def\ww{w(y)}
\def\ll{{\color{red}\ell}}
\def\ee{\mathrm{e}}
\def\startnewsection#1#2{ \section{#1}\label{#2}\setcounter{equation}{0}}   
\def\nnewpage{ }
\def\sgn{\mathop{\rm sgn\,}\nolimits}    
\def\Tr{\mathop{\rm Tr}\nolimits}    
\def\div{\mathop{\rm div}\nolimits}
\def\curl{\mathop{\rm curl}\nolimits}
\def\dist{\mathop{\rm dist}\nolimits}
\def\id{\mathop{\rm id}\nolimits}
\def\supp{\mathop{\rm supp}\nolimits}
\def\indeq{\quad{}}           
\def\period{.}                       
\def\semicolon{\,;}                  
\def\cmi#1{\text{~{{\coli IK: \underline{#1}}}~}}
\def\coli{\color{colorigor}}
\definecolor{colorigor}{rgb}{.5, 0.2, 0.8}
\def\colr{\color{red}}
\def\colrr{\color{black}}
\def\colb{\color{black}}
\def\coly{\color{lightgray}}
\definecolor{colorgggg}{rgb}{0.1,0.5,0.3}
\definecolor{colorllll}{rgb}{0.0,0.7,0.0}
\definecolor{colorhhhh}{rgb}{0.3,0.75,0.4}
\definecolor{colorpppp}{rgb}{0.7,0.0,0.2}
\definecolor{coloroooo}{rgb}{0.45,0.0,0.0}
\definecolor{colorqqqq}{rgb}{0.1,0.7,0}
\def\colg{\color{colorgggg}}
\def\collg{\color{colorllll}}
\def\cole{\color{coloroooo}}
\def\coleo{\color{colorpppp}}
\def\colu{\color{blue}}
\def\colc{\color{colorhhhh}}
\def\colW{\colb}   
\definecolor{coloraaaa}{rgb}{0.6,0.6,0.6}
\def\colw{\color{coloraaaa}}
\def\comma{ {\rm ,\qquad{}} }            
\def\commaone{ {\rm ,\quad{}} }          
\def\les{\lesssim}
\def\nts#1{{\color{red}\hbox{\bf ~#1~}}} 
\def\ntsf#1{\footnote{\color{colorgggg}\hbox{#1}}} 
\def\blackdot{{\color{red}{\hskip-.0truecm\rule[-1mm]{4mm}{4mm}\hskip.2truecm}}\hskip-.3truecm}
\def\bluedot{{\color{blue}{\hskip-.0truecm\rule[-1mm]{4mm}{4mm}\hskip.2truecm}}\hskip-.3truecm}
\def\purpledot{{\color{colorpppp}{\hskip-.0truecm\rule[-1mm]{4mm}{4mm}\hskip.2truecm}}\hskip-.3truecm}
\def\greendot{{\color{colorgggg}{\hskip-.0truecm\rule[-1mm]{4mm}{4mm}\hskip.2truecm}}\hskip-.3truecm}
\def\cyandot{{\color{cyan}{\hskip-.0truecm\rule[-1mm]{4mm}{4mm}\hskip.2truecm}}\hskip-.3truecm}
\def\reddot{{\color{red}{\hskip-.0truecm\rule[-1mm]{4mm}{4mm}\hskip.2truecm}}\hskip-.3truecm}
\def\tdot{{\color{green}{\hskip-.0truecm\rule[-.5mm]{3mm}{3mm}\hskip.2truecm}}\hskip-.1truecm}
\def\gdot{\greendot}
\def\bdot{\bluedot}
\def\ydot{\cyandot}
\def\rdot{\cyandot}
\def\fractext#1#2{{#1}/{#2}}
\def\ii{\hat\imath}
\def\fei#1{\textcolor{blue}{#1}}
\def\vlad#1{\textcolor{cyan}{#1}}
\def\igor#1{\text{{\textcolor{colorqqqq}{#1}}}}
\def\igorf#1{\footnote{\text{{\textcolor{colorqqqq}{#1}}}}}
\def\Omf{\Omega_{\text f}}
\def\Ome{\Omega_{\text e}}
\def\Omb{\Omega_{\text b}}
\def\Gaf{\Gamma_{\text f}}
\def\Gae{\Gamma_{\text e}}
\def\Gab{\Gamma_{\text b}}
\def\Gac{\Gamma_{\text c}}
\def\Nf{N^{\text f}}
\def\Ne{N^{\text e}}

\newcommand{\p}{\partial}
\renewcommand{\d}{\mathrm{d}}
\newcommand{\UE}{U^{\rm E}}
\newcommand{\PE}{P^{\rm E}}
\newcommand{\KP}{K_{\rm P}}
\newcommand{\uNS}{u^{\rm NS}}
\newcommand{\vNS}{v^{\rm NS}}
\newcommand{\pNS}{p^{\rm NS}}
\newcommand{\omegaNS}{\omega^{\rm NS}}
\newcommand{\uE}{u^{\rm E}}
\newcommand{\vE}{v^{\rm E}}
\newcommand{\pE}{p^{\rm E}}
\newcommand{\omegaE}{\omega^{\rm E}}
\newcommand{\ua}{u_{\rm   a}}
\newcommand{\va}{v_{\rm   a}}
\newcommand{\omegaa}{\omega_{\rm   a}}
\newcommand{\ue}{u_{\rm   e}}
\newcommand{\ve}{v_{\rm   e}}
\newcommand{\omegae}{\omega_{\rm e}}
\newcommand{\omegaeic}{\omega_{{\rm e}0}}
\newcommand{\ueic}{u_{{\rm   e}0}}
\newcommand{\veic}{v_{{\rm   e}0}}
\newcommand{\up}{u^{\rm P}}
\newcommand{\vp}{v^{\rm P}}
\newcommand{\tup}{{\tilde u}^{\rm P}}
\newcommand{\bvp}{{\bar v}^{\rm P}}
\newcommand{\omegap}{\omega^{\rm P}}
\newcommand{\tomegap}{\tilde \omega^{\rm P}}
\newcommand{\eps}{\varepsilon}  
\newcommand{\eqnb}{\begin{equation}}
\newcommand{\eqne}{\end{equation}}
  
\renewcommand{\up}{u^{\rm P}}
\renewcommand{\vp}{v^{\rm P}}
\renewcommand{\omegap}{\Omega^{\rm P}}
\renewcommand{\tomegap}{\omega^{\rm P}}

\begin{abstract}
We consider the Kuramoto-Sivashinsky equation
$u_t+uu_x+u_{xx}+u_{xxxx}=0$ with $L$-periodic boundary conditions and zero spatial mean.
We prove that every smooth solution satisfies $\|u(t)\|_{L^2}\leq CL^{5/6}+CL^{3/2}/t$ for $t>0$, where $C$ is a universal constant.
In particular, $\|u(t)\|_{L^2}\leq CL^{5/6}$ for all $t\geq L^{2/3}$, independently of the initial norm.
\end{abstract}

\maketitle

\setcounter{tocdepth}{2} 
\section{Introduction}\label{sec00}
We consider the Kuramoto-Sivashinsky (KS) equation
  \begin{equation}
   u_t
   + u u_x
   + u_{xx}
   + u_{xxxx}
   = 0
   \label{EQ01}
  \end{equation}
for a real-valued function $u=u(t,x)$, where $t\geq 0$ and $x$ ranges over $\TTL=\RR/L\ZZ$ with \mbox{$L>0$}.
Integrating~\eqref{EQ01} over one period shows that the spatial mean of $u$ is independent of time, and we assume that it vanishes, i.e.,
  \begin{equation}
   \int_0^L u(t,x)\,dx
   = 0
   \comma t\geq 0
   .
   \label{EQ02}
  \end{equation}
In this paper, we show that every smooth solution $u$ satisfies $\|u(t)\|_{L^2} \leq CL^{5/6} + CL^{3/2}/t$ for $t >  0$ and, moreover, that $\|u(t)\|_{L^2} \leq C L^{5/6}$ for all $t \geq L^{2/3}$.

The origins of equation~\eqref{EQ01} go back to Sivashinsky~\cite{S1}, who derived an asymptotic equation for the displacement of unstable flame fronts.
In one spatial dimension, neglecting thermal expansion and differentiating this equation yields~\eqref{EQ01} after rescaling. 
See also~\cite{S2} for a derivation under stoichiometric conditions. 
The equation was also derived independently by Kuramoto and Tsuzuki~\cite{KT} in the context of reaction--diffusion systems.
A central question in the analysis of the KS equation is how the energy $\|u(t)\|_{L^2}^2$ evolves in time depending on the spatial period $L$;
in particular, one may seek an absorbing ball in $L^2$ with a radius of the form $CL^\alpha$, where $C$ is independent of initial data and the period length $L$.
The conjectured optimal bound is $\limsup_{t \to \infty} \|u(t)\|_{L^2} \leq CL^{1/2}$, which would imply that the long-time energy per unit length remains bounded independently of domain size. 
The first long-time bound of the form $CL^\alpha$ was established by Nicolaenko, Scheurer, and Temam~\cite{NST}, who used a background-flow argument to achieve $C L^{5/2}$,
under the assumption that $u_0$ is odd.
Il'yashenko \cite{I} removed the oddness assumption by a different
argument based on the Burgers nonlinearity.
The bounds were later improved in~\cite{G},
\cite{CEES}, \cite{GO}, and~\cite{BG}, with the bounds
$C L^{5/2}$, $C L^{8/5}$, $o(L^{3/2})$, and $C L^{3/2}$, respectively.
Otto~\cite{O} subsequently obtained the long-time averaged estimate
\[
\limsup_{T\to\infty}
\left(\frac1T\int_0^T\|u(t)\|_{L^2(0,L)}^2\,dt\right)^{1/2}
\le C_\varepsilon L^{5/6}(\log (L+2))^{5/3+\varepsilon}
\]
for every $\varepsilon>0$, with $C_\varepsilon$ independent of $L$ and the initial data.
He also obtained bounds growing only as powers of $\log L$ for the normalized long-time averages of $\bigl|\,|\partial_x|^s u\bigr|^2$, for every $1/3<s\le2$.
The proof combines a priori estimates for the inhomogeneous inviscid Burgers equation with regularity estimates for the capillary Burgers equation, using the nonlinearity to obtain cubic control of spatial variations.
Goldman, Josien, and Otto~\cite{GJO} later simplified the Burgers estimate through a modified K\'arm\'an--Howarth--Monin identity and improved the logarithmic factor in the long-time averaged $L^2$ bound to $(\log L)^{1/2+\varepsilon}$.
For other results on the KS equation regarding backward behavior, steady states, well-posedness, analyticity, numerical studies of spatial scaling, and related topics, see~\cite{AM,BKRZ,KMal,KM,Mi,Mo,ST,T1,W,WH}.

In the present work, we remove the logarithmic loss and improve the long-time averaged $L^2$ bound to the pointwise-in-time estimate $\|u(t)\|_{L^2} \leq CL^{5/6}$ for all $t \geq L^{2/3}$, independently of initial data. 
The proof relies on regularization through the heat operator; we summarize the argument below.

In Section~\ref{sec01}, we state Theorem~\ref{T01}, which gives the $L^{5/6}$ bound. 
In Section~\ref{sec02}, we begin with an example which motivates the proof of Theorem~\ref{T01}. 
The main goal of the proof is to control the term $\|u_x\|_{L^2}^2$ in the energy identity~\eqref{EQ301}
with the cubic term $\|u_x\|_{L^3}^3$, which is obtained by differentiating~\eqref{EQ01} and testing with $|u_x|$.
However, this approach has two obstructions which prevent any useful estimates from closing. 
Instead, we apply $\partial_x P_s$ to~\eqref{EQ01}, where $P_s = e^{s\partial_x^2}$, and we average over the smoothing scale $s\in [0,S]$ with $S$ arbitrary. 
The resulting energy estimate is given in Proposition~\ref{P01}, whose proof is provided in Section~\ref{sec03}. 
Notably, this estimate admits a cubic dissipation $\frac{S}{4}\|\partial_x P_S u\|_{L^3}^3$ and, due to the properties of $P_s$, see Lemma~\ref{L00}, it includes a correction $H_S(u)$ to the energy which is bounded by $\frac14\|u\|_{L^2}^2$, so that the corrected energy is comparable to $\|u\|_{L^2}^2$.
Then, since $S$ is arbitrary, Lemma~\ref{L01} demonstrates that this cubic term controls the bad term $\|u_x\|_{L^2}^2$.
Finally, for the proof of Theorem~\ref{T01}, we average over intervals $[a,b]$ to obtain~\eqref{EQ322}, which shows that the energy at time $b$ depends sublinearly on its value at $a$.
Using~\eqref{EQ322} on suitable small intervals completes the proof of the theorem.
In Section~\ref{sec03}, we establish Proposition~\ref{P01}, the proof of which uses several convenient properties of the heat operator $P_s$ as well as a regularization $\Phi_\eps$ of the function $\Phi(z)=\frac12z|z|$ defining $H_S$ to justify computations for the transport term in~\eqref{EQ01}. 
The cubic term arises from the local variance $P_s(u^2)-(P_su)^2$, which measures the failure of $P_s$ to commute with the nonlinearity and is bounded from below by a reverse local Poincar\'e inequality, see Lemma~\ref{L02}.

\section{Main results}\label{sec01}
Functions on $\TTL$ are identified with $L$-periodic functions on~$\RR$.
We abbreviate $\|f\|=\|f\|_{L^2}$
and $(f,g)_L=\int_{0}^{L}f g\,dx$.
By $C$ we denote a positive universal constant whose value may change from line to line; it never depends on $L$, on the solution, or on~$t$.

The following theorem is the main result of the paper.

\begin{Theorem}
\label{T01}
Let $L>0$, and let $u$ be a smooth solution of~\eqref{EQ01} on $[0,\infty)\times\TTL$ satisfying~\eqref{EQ02}.
Then
  \begin{equation}
   \|u(t)\|_{L^2}
   \leq
   C L^{5/6}
   + \frac{CL^{3/2}}{t}
   \comma t>0
   .
   \label{EQ202}
  \end{equation}
Moreover, $\|u(t)\|_{L^2}^2\leq CL^{5/3}$ for every~$t\geq L^{2/3}$.
In particular,
  \begin{equation}
   \limsup_{t\to\infty}\|u(t)\|_{L^2}
   \leq
   C L^{5/6}
   .
   \label{EQ203}
  \end{equation}
\end{Theorem}


Theorem~\ref{T01} is proven in Sections~\ref{sec02} and~\ref{sec03}.

\section{Proof of Theorem~\ref{T01}}\label{sec02}
Throughout this section, $u$ denotes a smooth solution of~\eqref{EQ01} on~$[0,\infty)\times\TTL$ satisfying~\eqref{EQ02}.
Multiplying~\eqref{EQ01} by $u$, integrating over one period, and using $\int_0^Lu^2u_x\,dx=\frac13\int_0^L(u^3)_x\,dx=0$, we obtain the energy identity
  \begin{equation}
   \frac12\frac{d}{dt}\|u\|^2
   =
   \|u_x\|^2
   - \|u_{xx}\|^2
   .
   \label{EQ301}
  \end{equation}
The term $\|u_x\|^2$ is the obstruction to a bound on~$\|u\|$.
In Section~\ref{sec02a}, we explain in a motivating example how the nonlinearity can be used to control it.
In Section~\ref{sec02b}, this leads to a correction $H_S$ of the energy, defined through the heat operator, and to an energy inequality with a cubic dissipation term.
In Section~\ref{sec02d}, we use the cubic term to bound~$\|u\|$ and conclude the proof.

\subsection{A motivating example}\label{sec02a}
We want a cubic dissipation term, such as $\|u_x\|_{L^3}^3$, to control the term $\|u_x\|^2$ on the right side of~\eqref{EQ301}.
Such a term is produced by the nonlinearity at the level of the derivative.
Differentiating~\eqref{EQ01} in $x$, we obtain
$
   u_{xt}
   + (uu_x)_x
   + u_{xxx}
   + u_{xxxxx}
   = 0
$.
Testing this equation  with $|u_x|$, which is the derivative of $\Phi(z)=\frac12z|z|$ at $z=u_x$.
Since $(uu_x)_x=u_x^2+uu_{xx}$ and $|u_x|u_{xx}=\frac12\p_x(u_x|u_x|)$, an integration by parts gives
  \begin{equation}
   -\int_0^L|u_x|\,(uu_x)_x\,dx
   =
   -\int_0^L|u_x|^3\,dx
   + \frac12\int_0^Lu_x^2|u_x|\,dx
   =
   -\frac12\int_0^L|u_x|^3\,dx
   ,
   \label{EQ303}
  \end{equation}
and therefore
  \begin{equation}
   \frac{d}{dt}
   \Bigl(
    \frac12\int_0^Lu_x|u_x|\,dx
   \Bigr)
   + \frac12\int_0^L|u_x|^3\,dx
   =
   -\int_0^L|u_x|\,u_{xxx}\,dx
   - \int_0^L|u_x|\,u_{xxxxx}\,dx
   .
   \label{EQ304}
  \end{equation}
The cubic term on the left side of~\eqref{EQ304} has a favorable sign, while the first term on the right side equals $\int_0^L\sgn(u_x)u_{xx}^2\,dx$ and is thus bounded by~$\|u_{xx}\|^2$.
However, \eqref{EQ304} cannot be combined with~\eqref{EQ301}, for two reasons.
First, the functional $\frac12\int_0^Lu_x|u_x|\,dx$ is bounded by $\|u_x\|^2$ rather than $\|u\|^2$, so it cannot serve as a small correction to the energy.
Second, the last term in~\eqref{EQ304} contains five derivatives of~$u$; integrating it by parts requires $\|u_{xxx}\|^2$, which is not controlled by~\eqref{EQ301}, and produces a term concentrated on the set $\{u_x=0\}$ which has no sign.

\subsection{The energy inequality}\label{sec02b}
Both issues are removed by replacing $u_x$ with heat-smoothed derivatives and averaging over the smoothing scales.
For $s\geq 0$, denote by $P_s=\ee^{s\p_x^2}$ the heat operator on $\TTL$, i.e., the operator which multiplies the Fourier coefficient at the frequency $k\in(2\pi/L)\ZZ$ by~$\ee^{-sk^2}$.
The parameter $s$ is a smoothing scale and is unrelated to the time~$t$.
For $s>0$, we have $P_sf=p_s\ast f$, where the periodic heat kernel
  \begin{equation}
   p_s(x)
   =
   \frac{1}{\sqrt{4\pi s}}
   \sum_{n\in\ZZ}
   \ee^{-(x+nL)^2/(4s)}
   \llabel{EQ306}
  \end{equation}
is positive and satisfies $\int_0^L p_s(x)\,dx=1$.
We use the following properties of~$P_s$.

\begin{Lemma}
\label{L00}
Let $f\in C^\infty(\TTL)$ and $s,r\geq0$.
Then the following hold.
\begin{enumerate}[label=\textup{(\roman*)},ref=\textup{(\roman*)},leftmargin=*,itemsep=2pt]
\item\label{L00i}
$P_0f=f$, $P_sP_rf=P_{s+r}f$, and $\p_xP_sf=P_s\p_xf$.
\item\label{L00ii}
$\p_sP_sf=\p_x^2P_sf$.
\item\label{L00iii}
$\|P_sf\|\leq\|f\|$.
\item\label{L00iv}
If $u\in C^\infty([0,\infty)\times\TTL)$, then the function $(t,s,x)\mapsto P_su(t,\cdot)(x)$ belongs to $C^\infty([0,\infty)\times[0,\infty)\times\TTL)$.
\item\label{L00v}
For every $S>0$, we have the dissipation identity
  \begin{equation}
   \int_0^S\|\p_xP_sf\|^2\,ds
   =
   \frac12\|f\|^2
   - \frac12\|P_Sf\|^2
   \leq
   \frac12\|f\|^2
   .
   \label{EQ307}
  \end{equation}
\end{enumerate}
\end{Lemma}

\begin{proof}
Since $P_s$ multiplies the Fourier coefficient at $k$ by $\ee^{-sk^2}$, the properties~\ref{L00i} and~\ref{L00ii} follow from $\ee^{-sk^2}\ee^{-rk^2}=\ee^{-(s+r)k^2}$ and $\p_s\ee^{-sk^2}=(ik)^2\ee^{-sk^2}$, while~\ref{L00iii} follows from $\ee^{-sk^2}\leq1$ and Parseval's identity.
For~\ref{L00iv}, note that the Fourier coefficients of every $t$-derivative of $u$ decay faster than any power of $|k|$, uniformly for $t$ in compact sets; since $|\ee^{-sk^2}|\leq1$ and every derivative in $s$ or $x$ only contributes a power of $k$, the Fourier series of every derivative of $P_su$ in $(t,s,x)$ converges uniformly on compact subsets of $[0,\infty)\times[0,\infty)\times\TTL$.
Finally, by~\ref{L00ii} and an integration by parts, $\frac{d}{ds}\|P_sf\|^2=2\int_0^LP_sf\,\p_x^2P_sf\,dx=-2\|\p_xP_sf\|^2$, and integrating over $s\in[0,S]$ gives~\eqref{EQ307}.
\end{proof}

For $S>0$ and $f\in C^\infty(\TTL)$, we define
  \begin{equation}
   H_S(f)
   =
   \frac12
   \int_0^S\int_0^L
   \p_xP_sf\,|\p_xP_sf|
   \,dx\,ds
   \label{EQ308}
  \end{equation}
and we call $\|u\|^2+\frac12H_S(u)$ the \emph{corrected energy}.
By~\eqref{EQ307}, we have
  \begin{equation}
   |H_S(f)|
   \leq
   \frac12\int_0^S\|\p_xP_sf\|^2\,ds
   \leq
   \frac14\|f\|^2
   .
   \label{EQ309}
  \end{equation}
The corrected energy satisfies the following energy inequality, which replaces the combination of~\eqref{EQ301} and~\eqref{EQ304}.

\begin{Proposition}
\label{P01}
Let $u$ be a smooth solution of~\eqref{EQ01} on $[0,\infty)\times\TTL$, and let $S>0$.
Then the function $t\mapsto\|u(t)\|^2+\frac12H_S(u(t))$ belongs to $C^1([0,\infty))$ and satisfies
  \begin{equation}
   \frac78\|u\|^2
   \leq
   \|u\|^2
   + \frac12 H_S(u)
   \leq
   \frac98\|u\|^2
   \label{EQ310}
  \end{equation}
and
  \begin{equation}
   \frac{d}{dt}
   \Bigl(
    \|u\|^2
    + \frac12 H_S(u)
   \Bigr)
   + \|u_{xx}\|^2
   + \frac{S}{4}\|\p_xP_Su\|_{L^3}^3
   \leq
   \frac52\|u_x\|^2
   \comma t\geq0
   .
   \label{EQ311}
  \end{equation}
\end{Proposition}

Compared with~\eqref{EQ304}, the cubic term now involves the smoothed gradient at the scale~$S$, and the right side contains only~$\|u_x\|^2$.
We postpone the proof of Proposition~\ref{P01} to Section~\ref{sec03}.

\subsection{Conclusion of the proof}\label{sec02d}
The cubic term in~\eqref{EQ311} appears at a single scale $S$, but $S>0$ is arbitrary.
We use ordinary time averages.
For an interval $[a,b]$, write
  \begin{equation}
   \langle F\rangle
   =
   \frac{1}{b-a}\int_a^b F(t)\,dt
   .
   \llabel{EQ312}
  \end{equation}
The next lemma shows that a bound on the averaged cubic term, uniform in the smoothing scale, controls the troublesome term $\langle\|u_x\|^2\rangle$ as well as the averaged energy.

\begin{Lemma}
\label{L01}
Let $0\leq a<b$, and let $u\in C^\infty([a,b]\times\TTL)$ have zero spatial mean.
Define $K\in[0,\infty]$ by
  \begin{equation}
   K^3
   =
   \sup_{s>0}\,
   s\,\langle\|\p_xP_su\|_{L^3}^3\rangle
   .
   \label{EQ313}
  \end{equation}
Then, for every $\delta>0$,
  \begin{equation}
   \langle\|u_x\|^2\rangle
   \leq
   \delta\langle\|u_{xx}\|^2\rangle
   + \ee^2\delta^{-2/3}L^{1/3}K^2
   ,
   \label{EQ314}
  \end{equation}
and
  \begin{equation}
   \langle\|u\|^2\rangle
   \leq
   7LK^2
   .
   \label{EQ315}
  \end{equation}
\end{Lemma}

\begin{proof}
Let $s>0$.
H\"older's inequality and~\eqref{EQ313} give
  \begin{equation}
   \langle\|\p_xP_su\|^2\rangle
   \leq
   L^{1/3}\langle\|\p_xP_su\|_{L^3}^3\rangle^{2/3}
   \leq
   L^{1/3}K^2s^{-2/3}
   .
   \label{EQ316}
  \end{equation}
To prove~\eqref{EQ314}, note that
  \begin{equation}
   k^2
   \leq
   \delta k^4
   + k^2\ee^{2-2\delta k^2}
   \comma k\in\RR
   .
   \llabel{EQ317}
  \end{equation}
Indeed, if $\delta k^2\geq1$, then $k^2\leq\delta k^4$, while if $\delta k^2<1$, then $\ee^{2-2\delta k^2}>1$.
Multiplying by $|\widehat u(t,k)|^2$, summing over $k$, and using Parseval's identity, we obtain
  \begin{equation}
   \|u_x\|^2
   \leq
   \delta\|u_{xx}\|^2
   + \ee^2\|\p_xP_\delta u\|^2
   .
   \llabel{EQ318}
  \end{equation}
Average over $[a,b]$ and use~\eqref{EQ316} with~$s=\delta$.

To prove~\eqref{EQ315}, note that since $u$ has zero mean, only frequencies with $|k|\geq2\pi/L$ occur, and thus $\|P_{L^2}u\|^2\leq\ee^{-8\pi^2}\|u\|^2$.
By~\eqref{EQ307},
  \begin{equation}
   \frac{1-\ee^{-8\pi^2}}{2}\|u\|^2
   \leq
   \int_0^{L^2}\|\p_xP_su\|^2\,ds
   .
   \llabel{EQ319}
  \end{equation}
Averaging and using~\eqref{EQ316}, we get
  \begin{equation}
   \frac{1-\ee^{-8\pi^2}}{2}\langle\|u\|^2\rangle
   \leq
   L^{1/3}K^2\int_0^{L^2}s^{-2/3}\,ds
   =
   3LK^2
   .
   \llabel{EQ320}
  \end{equation}
Since $6/(1-\ee^{-8\pi^2})<7$, this proves~\eqref{EQ315}.
\end{proof}

To prove Theorem~\ref{T01}, we first integrate~\eqref{EQ311} over a time interval.
The cubic term bounds the average energy, which gives a time at which the energy is small, and the same integrated inequality bounds the growth of the energy from that time to the end of the interval.
In effect, we obtain~\eqref{EQ322}, which states that the energy at the end of an interval depends only sublinearly on its value at the beginning. 
We then apply this bound on suitably chosen shorter intervals to remove the dependence on the initial energy.

\begin{proof}[Proof of Theorem~\ref{T01}]
Write
  \begin{equation}
   E(t)
   =
   \|u(t)\|^2
   .
   \llabel{EQ321}
  \end{equation}
We first claim that, for $0\leq a<b$ and with $T=b-a$,
  \begin{equation}
   E(b)
   \leq
   CL\left(\frac{E(a)}{T}\right)^{2/3}
   + CL^{5/3}
   + CLT
   .
   \label{EQ322}
  \end{equation}

To prove~\eqref{EQ322}, we take $\langle\cdot\rangle$ and $K$ over the interval~$[a,b]$.
Integrating~\eqref{EQ311} and using~\eqref{EQ310}, we obtain, for every $S>0$,
  \begin{equation}
   \frac78E(b)
   + \int_a^b\|u_{xx}\|^2\,dt
   + \frac{S}{4}\int_a^b\|\p_xP_Su\|_{L^3}^3\,dt
   \leq
   \frac98E(a)
   + \frac52\int_a^b\|u_x\|^2\,dt
   .
   \llabel{EQ323}
  \end{equation}
The right side is finite and does not depend on~$S$.
Taking the supremum over $S>0$ shows that $K<\infty$ and gives
  \begin{equation}
   \frac78E(b)
   + T\langle\|u_{xx}\|^2\rangle
   + \frac{T}{4}K^3
   \leq
   \frac98E(a)
   + \frac52T\langle\|u_x\|^2\rangle
   .
   \label{EQ324}
  \end{equation}
By~\eqref{EQ314} with $\delta=1/5$,
  \begin{equation}
   \frac52T\langle\|u_x\|^2\rangle
   \leq
   \frac{T}{2}\langle\|u_{xx}\|^2\rangle
   + CTL^{1/3}K^2
   .
   \llabel{EQ325}
  \end{equation}
Insert this into~\eqref{EQ324} and use
  \begin{equation}
   CL^{1/3}K^2
   \leq
   \frac18K^3
   + CL
   .
   \llabel{EQ326}
  \end{equation}
After multiplication by $8$, we obtain
  \begin{equation}
   7E(b)
   + 4\int_a^b\|u_{xx}\|^2\,dt
   + TK^3
   \leq
   9E(a)
   + CLT
   .
   \llabel{EQ327}
  \end{equation}
In particular,
  \begin{equation}
   E(b)
   \leq
   \frac97E(a)
   + CL(b-a)
   ,
   \label{EQ328}
  \end{equation}
and
  \begin{equation}
   K^3
   \leq
   \frac{9E(a)}{T}
   + CL
   .
   \label{EQ329}
  \end{equation}
The constants in these estimates are independent of $a$, $b$, and~$T$.
By~\eqref{EQ315} and~\eqref{EQ329},
  \begin{equation}
   \frac{1}{T}\int_a^bE(t)\,dt
   \leq
   7LK^2
   \leq
   CL\left(\frac{E(a)}{T}\right)^{2/3}
   + CL^{5/3}
   .
   \llabel{EQ330}
  \end{equation}
Since $E$ is continuous, there exists $t_0\in[a,b]$ at which $E$ is bounded by its average, and therefore
  \begin{equation}
   E(t_0)
   \leq
   CL\left(\frac{E(a)}{T}\right)^{2/3}
   + CL^{5/3}
   .
   \llabel{EQ331}
  \end{equation}
Applying~\eqref{EQ328} on $[t_0,b]$, which is immediate if $t_0=b$, we obtain
  \begin{equation}
   E(b)
   \leq
   \frac97E(t_0)
   + CL(b-t_0)
   .
   \llabel{EQ332}
  \end{equation}
Inserting the above bound on $E(t_0)$ and using $b-t_0\leq T$, we obtain~\eqref{EQ322}.

Next, we claim that, for $a\geq0$ and $0<T\leq L^{2/3}$,
  \begin{equation}
   E(a+T)
   \leq
   \frac{CL^3}{T^2}
   + CL^{5/3}
   .
   \label{EQ333}
  \end{equation}
To prove~\eqref{EQ333}, let
  \begin{equation}
   M
   =
   \max_{0\leq r\leq T}r^2E(a+r)
   ,
   \llabel{EQ334}
  \end{equation}
which is finite because $E$ is continuous.
If $M=0$, then $E(a+T)=0$ and there is nothing to prove.
Otherwise, the maximum is attained at some~$\sigma\in(0,T]$, and
  \begin{equation}
   M
   =
   \sigma^2E(a+\sigma)
   \comma
   \left(\frac{\sigma}{2}\right)^2E\left(a+\frac{\sigma}{2}\right)
   \leq
   M
   .
   \llabel{EQ335}
  \end{equation}
It follows that
  \begin{equation}
   E\left(a+\frac{\sigma}{2}\right)
   \leq
   4E(a+\sigma)
   .
   \llabel{EQ336}
  \end{equation}
We apply~\eqref{EQ322} on the interval $[a+\sigma/2,a+\sigma]$, whose length is~$\sigma/2$, and obtain
  \begin{align}
   \begin{split}
    E(a+\sigma)
    &\leq
    CL\left(\frac{E(a+\sigma/2)}{\sigma/2}\right)^{2/3}
    + CL^{5/3}
    + CL\sigma
    \\&\leq
    CL\sigma^{-2/3}E(a+\sigma)^{2/3}
    + CL^{5/3}
    ,
   \end{split}
   \llabel{EQ337}
  \end{align}
where in the last step we used $\sigma\leq T\leq L^{2/3}$.
Multiplying by $\sigma^2$ and using $M=\sigma^2E(a+\sigma)$ and Young's inequality, we get
  \begin{equation}
   M
   \leq
   CLM^{2/3}
   + CL^{5/3}T^2
   \leq
   \frac12M
   + CL^3
   + CL^{5/3}T^2
   ,
   \llabel{EQ338}
  \end{equation}
and hence $M\leq CL^3+CL^{5/3}T^2$.
Since $T^2E(a+T)\leq M$, the proof of the claim~\eqref{EQ333} is complete.

Finally, for $T\leq L^{2/3}$, we have $L^{5/3}\leq L^3/T^2$, and thus~\eqref{EQ333} gives
  \begin{equation}
   E(a+T)
   \leq
   \frac{CL^3}{T^2}
   \comma a\geq0
   \commaone 0<T\leq L^{2/3}
   .
   \llabel{EQ339}
  \end{equation}
If $0<t\leq L^{2/3}$, we take $a=0$ and $T=t$.
If $t\geq L^{2/3}$, we take $a=t-L^{2/3}$ and~$T=L^{2/3}$.
We obtain
  \begin{equation}
   E(t)
   \leq
   \begin{cases}
    CL^3/t^2
    ,
    & 0<t\leq L^{2/3}
    ,
    \\
    CL^{5/3}
    ,
    & t\geq L^{2/3}
    ,
   \end{cases}
   \llabel{EQ340}
  \end{equation}
which proves~\eqref{EQ202}, the bound for every $t\geq L^{2/3}$, and~\eqref{EQ203}.
\end{proof}

\begin{Remark}
\label{R03}
{\rm
The same time-interval argument also gives
  \begin{equation}
   \|u(t)\|^2
   \leq
   CL^{5/3}
   + C\ee^{-t/L^{2/3}}\|u(0)\|^2
   \comma t\geq0
   .
   \label{EQ341}
  \end{equation}
Indeed, put $T_0=L^{2/3}$ in~\eqref{EQ322}.
Young's inequality, with its coefficient chosen as $\ee^{-1}$, gives
  \begin{equation}
   E(a+T_0)
   \leq
   CL^{5/9}E(a)^{2/3}
   + CL^{5/3}
   \leq
   \ee^{-1}E(a)
   + CL^{5/3}
   .
   \llabel{EQ342}
  \end{equation}
Iteration at the times $nT_0$ yields
  \begin{equation}
   E(nT_0)
   \leq
   \ee^{-n}E(0)
   + CL^{5/3}
   .
   \llabel{EQ343}
  \end{equation}
For $nT_0\leq t<(n+1)T_0$, apply~\eqref{EQ328} to $[nT_0,t]$ and use $\ee^{-n}\leq\ee\,\ee^{-t/T_0}$ to obtain~\eqref{EQ341}.
The improvement in~\eqref{EQ202} is the absence of the initial norm: it gives the uniform bound $E(t)\leq CL^{5/3}$ for every $t\geq L^{2/3}$, even when $E(0)$ is arbitrarily large.
}
\end{Remark}

\section{Proof of Proposition~\ref{P01}}\label{sec03}
Throughout this section, $u$ is a smooth solution of~\eqref{EQ01} on $[0,\infty)\times\TTL$, the scale $S>0$ is fixed, and
  \begin{equation}
   w(t,s,x)
   =
   \p_xP_su(t,\cdot)(x)
   \comma (t,s,x)\in[0,\infty)\times[0,\infty)\times\TTL
   ,
   \llabel{EQ401}
  \end{equation}
so that $w(t,0,\cdot)=u_x(t,\cdot)$.
Thus $H_S(u)=\int_0^S\int_0^L\Phi(w)\,dx\,ds$ with $\Phi(z)=\frac12z|z|$, i.e., $H_S(u)$ is obtained by integrating the functional in~\eqref{EQ304}, with $u_x$ replaced by $w$, over the scales $s\in[0,S]$.
The averaging in $s$ gains one derivative: by~\eqref{EQ309}, $H_S(u)$ is bounded by a fraction of $\|u\|^2$, which resolves the first issue raised in Section~\ref{sec02a} and gives the bounds~\eqref{EQ310}.
Since adding one half of
  \begin{equation}
   \frac{d}{dt}H_S(u)
   + \frac{S}{2}\|\p_xP_Su\|_{L^3}^3
   \leq
   \|u_x\|^2
   + 2\|u_{xx}\|^2
   \label{EQ402}
  \end{equation}
to the energy identity~\eqref{EQ301} gives~\eqref{EQ311}, it remains to show that $t\mapsto H_S(u(t))$ belongs to $C^1([0,\infty))$ and that~\eqref{EQ402}, the smoothed counterpart of~\eqref{EQ304}, holds.

\subsection{The smoothed equation}\label{sec03a}
Applying $\p_xP_s$ to~\eqref{EQ01} and using $\p_sw=w_{xx}$ (see Lemma~\ref{L00}\ref{L00i}--\ref{L00ii}), which implies $w_{xxxx}=\p_s^2w$, we obtain
  \begin{equation}
   w_t
   + \p_xP_s(uu_x)
   + \p_sw
   + \p_s^2w
   = 0
   \inon{on $[0,\infty)\times[0,\infty)\times\TTL$.}
   \label{EQ403}
  \end{equation}
In~\eqref{EQ403}, the two linear terms are derivatives in the smoothing variable, so that after an integration in $s$ they only contribute boundary terms at $s=0$ and $s=S$; this resolves the second issue raised in Section~\ref{sec02a}.
However, we pay the price in the transport term.
Since $P_s$ does not commute with multiplication, $\p_xP_s(uu_x)=\frac12\p_x^2P_s(u^2)$ differs from $\p_x(P_su\,w)=\frac12\p_x^2((P_su)^2)$ by the term $\frac12\p_x^2\tau$, where $\tau=P_s(u^2)-(P_su)^2$ is the local variance of~$u$.
We need two properties of the local variance, collected in the next lemma.
The second one is the reverse local Poincar\'e inequality, which bounds the local variance at the scale $\sqrt S$ from below by the square of the smoothed gradient.
Its proof uses that, since $p_s$ is a probability density, the Cauchy-Schwarz inequality gives
  \begin{equation}
   (P_sf)^2
   \leq
   P_s(f^2)
   \inon{on $\TTL$.}
   \label{EQ404}
  \end{equation}

\begin{Lemma}
\label{L02}
Let $f\in C^\infty(\TTL)$.
Then the local variance
  \begin{equation}
   \tau(s,x)
   =
   P_s(f^2)(x)
   - \bigl(P_sf(x)\bigr)^2
   \comma (s,x)\in[0,\infty)\times\TTL
   ,
   \llabel{EQ405}
  \end{equation}
satisfies
  \begin{equation}
   \p_s\tau
   =
   \p_x^2\tau
   + 2(\p_xP_sf)^2
   \comma
   \tau(0,\cdot)
   = 0
   ,
   \label{EQ406}
  \end{equation}
and, for every $S>0$,
  \begin{equation}
   \tau(S,x)
   \geq
   2S\bigl(\p_xP_Sf(x)\bigr)^2
   \comma x\in\TTL
   .
   \label{EQ407}
  \end{equation}
In particular, $\tau\geq 0$.
\end{Lemma}

\begin{proof}
We have $\p_sP_s(f^2)=\p_x^2P_s(f^2)$ and
  \begin{equation}
   \p_s(P_sf)^2
   =
   2P_sf\,\p_x^2P_sf
   =
   \p_x^2\bigl((P_sf)^2\bigr)
   - 2(\p_xP_sf)^2
   .
   \llabel{EQ408}
  \end{equation}
Subtracting the two identities gives the equation in~\eqref{EQ406}, while $\tau(0,\cdot)=0$ since $P_0f=f$.
To prove~\eqref{EQ407}, we note that, by~\eqref{EQ406},
  \begin{equation}
   \p_r\bigl(P_{S-r}\tau(r,\cdot)\bigr)
   =
   -\p_x^2P_{S-r}\tau(r,\cdot)
   + P_{S-r}\p_r\tau(r,\cdot)
   =
   2P_{S-r}\bigl((\p_xP_rf)^2\bigr)
   \comma r\in[0,S]
   .
   \llabel{EQ409}
  \end{equation}
Integrating over $r\in[0,S]$ and using $\tau(0,\cdot)=0$, we obtain
  \begin{equation}
   \tau(S,\cdot)
   =
   2\int_0^S P_{S-r}\bigl((\p_xP_rf)^2\bigr)\,dr
   .
   \label{EQ410}
  \end{equation}
By~\eqref{EQ404} and the semigroup property, $P_{S-r}((\p_xP_rf)^2)\geq(P_{S-r}\p_xP_rf)^2=(\p_xP_Sf)^2$ for every $r\in[0,S]$, and inserting this into~\eqref{EQ410} gives~\eqref{EQ407}.
\end{proof}

\subsection{Regularization}\label{sec03b}
Since $\Phi(z)=\frac12z|z|$ belongs to $C^1(\RR)$ with $\Phi'(z)=|z|$, and since $w$ is smooth in $(t,s,x)$ by Lemma~\ref{L00}\ref{L00iv}, differentiation under the integral sign shows that $t\mapsto H_S(u(t))$ belongs to $C^1([0,\infty))$, with
  \begin{equation}
   \frac{d}{dt}H_S(u)
   =
   \int_0^S\int_0^L|w|\,w_t\,dx\,ds
   .
   \label{EQ411}
  \end{equation}
Thus~\eqref{EQ402} should follow by testing~\eqref{EQ403} with $|w|$ and integrating over $[0,S]\times\TTL$, just as~\eqref{EQ304} was obtained by testing
the differentiated Kuramoto-Sivashinsky equation 
with~$|u_x|$.
Here we meet a regularity issue: the computation of the transport term involves two derivatives of the multiplier, and $|z|$ is not differentiable at the origin.
The fix is to approximate $|z|$ by the smooth convex function $\sqrt{\eps^2+z^2}-\eps$, where $\eps>0$, i.e., to replace $\Phi$ by
  \begin{equation}
   \Phi_\eps(z)
   =
   \int_0^z
   \Bigl(
    \sqrt{\eps^2+y^2}
    - \eps
   \Bigr)
   \,dy
   \comma z\in\RR
   .
   \llabel{EQ412}
  \end{equation}
Then $\Phi_\eps\in C^\infty(\RR)$, $\Phi_\eps(0)=0$, and
  \begin{equation}
   0
   \leq
   \Phi_\eps'(z)
   \leq
   |z|
   \commaone
   \bigl|\Phi_\eps'(z)-|z|\bigr|
   \leq
   \eps
   \commaone
   |\Phi_\eps''(z)|
   =
   \frac{|z|}{\sqrt{\eps^2+z^2}}
   \leq
   1
   \commaone
   \Phi_\eps'''(z)
   =
   \frac{\eps^2}{(\eps^2+z^2)^{3/2}}
   >
   0
   .
   \label{EQ413}
  \end{equation}
Since $\Phi_\eps'$ is even, $\Phi_\eps$ is odd, and the first inequality in~\eqref{EQ413} gives~$|\Phi_\eps(z)|\leq\frac12z^2$.
The next lemma is the energy inequality for~$\Phi_\eps$.

\begin{Lemma}
\label{L03}
Let $\eps>0$, and let $\tau$ denote the local variance of $u(t,\cdot)$.
Then the function $t\mapsto\int_0^S\int_0^L\Phi_\eps(w)\,dx\,ds$ belongs to $C^1([0,\infty))$ and satisfies
  \begin{equation}
   \frac{d}{dt}\int_0^S\int_0^L\Phi_\eps(w)\,dx\,ds
   \leq
   -\frac14\int_0^L\tau(S,x)\,\Phi_\eps'\bigl(\p_xP_Su(x)\bigr)\,dx
   + \|u_x\|^2
   + 2\|u_{xx}\|^2
   \comma t\geq0
   .
   \label{EQ414}
  \end{equation}
\end{Lemma}

\begin{proof}
\textit{Step~1. Testing the equation.}
As for $H_S$, the function $t\mapsto\int_0^S\int_0^L\Phi_\eps(w)\,dx\,ds$ is $C^1$, and testing~\eqref{EQ403} with $\Phi_\eps'(w)$ gives
  \begin{align}
   \begin{split}
    \frac{d}{dt}\int_0^S\int_0^L\Phi_\eps(w)\,dx\,ds
    &=
    -\int_0^S\int_0^L\Phi_\eps'(w)\,\p_xP_s(uu_x)\,dx\,ds
    - \int_0^S\int_0^L\Phi_\eps'(w)\,\p_sw\,dx\,ds
    \\&\indeq
    - \int_0^S\int_0^L\Phi_\eps'(w)\,\p_s^2w\,dx\,ds
    =
    I_1
    + I_2
    + I_3
    .
   \end{split}
   \label{EQ415}
  \end{align}

\textit{Step~2. The linear terms.}
Since $\Phi_\eps'(w)\p_sw=\p_s(\Phi_\eps(w))$ and $w(t,0,\cdot)=u_x$, we have
  \begin{equation}
   I_2
   =
   -\int_0^S\p_s\int_0^L\Phi_\eps(w)\,dx\,ds
   =
   \int_0^L\Phi_\eps(u_x)\,dx
   - \int_0^L\Phi_\eps(\p_xP_Su)\,dx
   .
   \llabel{EQ416}
  \end{equation}
By $|\Phi_\eps(z)|\leq\frac12z^2$ and $\|\p_xP_Su\|=\|P_Su_x\|\leq\|u_x\|$, each of the two integrals is bounded in absolute value by $\frac12\|u_x\|^2$, and thus
  \begin{equation}
   I_2
   \leq
   \|u_x\|^2
   .
   \label{EQ417}
  \end{equation}
For $I_3$, we integrate by parts in $s$ and use $\p_sw=w_{xx}$, obtaining
  \begin{equation}
   I_3
   =
   -\int_0^L\Phi_\eps'(w)\,w_{xx}\,dx\,\bigg|_{s=0}^{s=S}
   + \int_0^S\int_0^L\Phi_\eps''(w)\,w_{xx}^2\,dx\,ds
   .
   \llabel{EQ418}
  \end{equation}
Using $\int_0^L\Phi_\eps'(w)w_{xx}\,dx=-\int_0^L\Phi_\eps''(w)w_x^2\,dx$, together with $w_x(t,0,\cdot)=u_{xx}$, $w_x(t,S,\cdot)=P_Su_{xx}$, and $w_{xx}=\p_xP_su_{xx}$, we get
  \begin{align}
   \begin{split}
    I_3
    &=
    -\int_0^L\Phi_\eps''(u_x)\,u_{xx}^2\,dx
    + \int_0^L\Phi_\eps''(\p_xP_Su)\,(P_Su_{xx})^2\,dx
    \\&\indeq
    + \int_0^S\int_0^L\Phi_\eps''(w)\,(\p_xP_su_{xx})^2\,dx\,ds
    .
   \end{split}
   \llabel{EQ419}
  \end{align}
By $|\Phi_\eps''|\leq 1$, \eqref{EQ307} with $u_{xx}$ in place of $f$, and $\|P_Su_{xx}\|\leq\|u_{xx}\|$, we conclude that
  \begin{equation}
   I_3
   \leq
   \|u_{xx}\|^2
   + \|P_Su_{xx}\|^2
   + \frac12\|u_{xx}\|^2
   - \frac12\|P_Su_{xx}\|^2
   \leq
   2\|u_{xx}\|^2
   .
   \label{EQ420}
  \end{equation}

\textit{Step~3. The transport term.}
Fix $s\in[0,S]$; we suppress the arguments $(t,s,x)$.
Since $P_s(u^2)=(P_su)^2+\tau$ and $\p_xP_su=w$, we have
  \begin{equation}
   \p_xP_s(uu_x)
   =
   \frac12\p_x^2P_s(u^2)
   =
   \frac12\p_x^2\bigl((P_su)^2\bigr)
   + \frac12\tau_{xx}
   =
   w^2
   + (P_su)\,w_x
   + \frac12\tau_{xx}
   .
   \llabel{EQ421}
  \end{equation}
We multiply this identity by $-\Phi_\eps'(w)$ and integrate in~$x$.
Since $\Phi_\eps'(w)w_x=\p_x(\Phi_\eps(w))$, one integration by parts in the second term and two in the third term give
  \begin{align}
   \begin{split}
    -\int_0^L\Phi_\eps'(w)\,\p_xP_s(uu_x)\,dx
    &=
    \int_0^L
    \bigl(
     w\Phi_\eps(w)
     - w^2\Phi_\eps'(w)
    \bigr)
    \,dx
    - \frac12\int_0^L\tau\,\Phi_\eps'''(w)\,w_x^2\,dx
    \\&\indeq
    - \frac12\int_0^L\tau\,\Phi_\eps''(w)\,w_{xx}\,dx
    .
   \end{split}
   \label{EQ422}
  \end{align}
When $\tau=0$ and $\Phi_\eps(z)$ is replaced by $\frac12z|z|$, the first integral on the right side equals $-\frac12\int_0^L|w|^3\,dx$, as in~\eqref{EQ303}; the remaining two terms come from the commutator.
The integrand in the last term has no sign, and we write it as a derivative in the smoothing variable.
Since $\p_sw=w_{xx}$ and since $\tau$ satisfies~\eqref{EQ406}, two integrations by parts in the term containing $\tau_{xx}$ give
  \begin{align}
   \begin{split}
    \p_s\int_0^L\tau\,\Phi_\eps'(w)\,dx
    &=
    \int_0^L
    \Bigl(
     \bigl(\tau_{xx}+2w^2\bigr)\Phi_\eps'(w)
     + \tau\,\Phi_\eps''(w)\,w_{xx}
    \Bigr)
    \,dx
    \\&=
    \int_0^L
    \Bigl(
     \tau\,\Phi_\eps'''(w)\,w_x^2
     + 2\tau\,\Phi_\eps''(w)\,w_{xx}
     + 2w^2\Phi_\eps'(w)
    \Bigr)
    \,dx
    .
   \end{split}
   \llabel{EQ423}
  \end{align}
Solving for the integral of $\tau\Phi_\eps''(w)w_{xx}$ and substituting the result into~\eqref{EQ422}, we obtain
  \begin{align}
   \begin{split}
    -\int_0^L\Phi_\eps'(w)\,\p_xP_s(uu_x)\,dx
    &=
    -\int_0^L
    \Bigl(
     \frac12 w^2\Phi_\eps'(w)
     - w\Phi_\eps(w)
    \Bigr)
    \,dx
    - \frac14\int_0^L\tau\,\Phi_\eps'''(w)\,w_x^2\,dx
    \\&\indeq
    - \frac14\,\p_s\int_0^L\tau\,\Phi_\eps'(w)\,dx
    .
   \end{split}
   \llabel{EQ424}
  \end{align}
Integrating over $s\in[0,S]$ and using $\tau(0,\cdot)=0$, we arrive at
  \begin{align}
   \begin{split}
    I_1
    &=
    -\int_0^S\int_0^L
    \Bigl(
     \frac12 w^2\Phi_\eps'(w)
     - w\Phi_\eps(w)
    \Bigr)
    \,dx\,ds
    - \frac14\int_0^S\int_0^L
    \tau\,\Phi_\eps'''(w)\,w_x^2
    \,dx\,ds
    \\&\indeq
    - \frac14\int_0^L
    \tau(S,x)\,\Phi_\eps'\bigl(w(t,S,x)\bigr)
    \,dx
    .
   \end{split}
   \label{EQ425}
  \end{align}

\textit{Step~4. The signs in the transport term.}
Since $\tau\geq0$ by Lemma~\ref{L02} and $\Phi_\eps'''>0$, the second term on the right side of~\eqref{EQ425} is nonpositive.
For the first term, we set $G(z)=z\Phi_\eps'(z)-2\Phi_\eps(z)$.
Then $G(0)=0$ and
  \begin{equation}
   G'(z)
   =
   z\Phi_\eps''(z)
   - \Phi_\eps'(z)
   =
   \eps
   - \frac{\eps^2}{\sqrt{\eps^2+z^2}}
   \geq
   0
   \comma z\in\RR
   ,
   \llabel{EQ426}
  \end{equation}
so that $G(z)$ has the sign of $z$, and $\frac12z^2\Phi_\eps'(z)-z\Phi_\eps(z)=\frac12zG(z)\geq0$.
Hence the first term on the right side of~\eqref{EQ425} is also nonpositive.
Combining~\eqref{EQ415}, \eqref{EQ417}, \eqref{EQ420}, and~\eqref{EQ425}, we obtain~\eqref{EQ414}.
\end{proof}

\subsection{The approximation argument}\label{sec03c}
We now let $\eps\to0$ in~\eqref{EQ414} and conclude the proof of Proposition~\ref{P01}.

\begin{proof}[Proof of Proposition~\ref{P01}]
As noted at the beginning of this section, it suffices to prove that $t\mapsto H_S(u(t))$ is $C^1$, which was shown in~\eqref{EQ411}, and that~\eqref{EQ402} holds.
We fix $t$ and let $\eps\to0$.
By~\eqref{EQ411}, the analogous formula for $\Phi_\eps$, and $|\Phi_\eps'(z)-|z||\leq\eps$, we have
  \begin{equation}
   \biggl|
    \frac{d}{dt}\int_0^S\int_0^L\Phi_\eps(w)\,dx\,ds
    - \frac{d}{dt}H_S(u)
   \biggr|
   =
   \biggl|
    \int_0^S\int_0^L
    \bigl(
     \Phi_\eps'(w)
     - |w|
    \bigr)
    \,w_t
    \,dx\,ds
   \biggr|
   \leq
   \eps\int_0^S\|w_t\|_{L^1}\,ds
   \to
   0
   ,
   \llabel{EQ427}
  \end{equation}
while $\Phi_\eps'(\p_xP_Su)\to|\p_xP_Su|$ uniformly on~$\TTL$.
Passing to the limit in~\eqref{EQ414} and then using~\eqref{EQ407} in the form $\tau(S,x)\,|\p_xP_Su(x)|\geq 2S|\p_xP_Su(x)|^3$, we obtain~\eqref{EQ402}.
\end{proof}

\begin{Remark}
\label{R04}
{\rm
For the limiting function $\Phi(z)=\frac12z|z|$, the expression $\frac12z^2\Phi'(z)-z\Phi(z)$ vanishes identically.
Thus, in contrast with~\eqref{EQ303}, the cubic term in~\eqref{EQ402} does not come from the first term on the right side of~\eqref{EQ425}.
It comes from the boundary term at $s=S$, which is generated by the commutator~$\tau$.
}
\end{Remark}


\par\medskip\noindent\textbf{Acknowledgments.} The authors were supported in part by the NSF grant DMS-2205493
and the Simons grant SFI-MPS-TSM-00025796.
The authors used ChatGPT~6.1 and Claude Opus~5.5 for help with mathematical arguments and writing.

%
\end{document}